\documentclass[10pt,leqno]{article}
\usepackage{latexsym}
\usepackage{amsthm}
\usepackage{amsfonts}
\usepackage{amssymb}
\usepackage{amsmath} 
\usepackage{graphicx}
\usepackage{url} 
\usepackage{color}
\usepackage{indentfirst}

\usepackage[top=30truemm,bottom=30truemm,left=30truemm,right=30truemm]{geometry}

\theoremstyle{plain}
  \newtheorem{theorem}{Theorem}[section]

  \newtheorem{lemma}{Lemma}[section]

\theoremstyle{remark}
  \newtheorem{remark}{Remark}[section]

\theoremstyle{definition}

\makeatletter

\@addtoreset{equation}{section}
\makeatother

\begin{document}

\markboth{Hideo Takaoka}{}

\title{Bilinear Strichartz estimates for the KdV equation on the torus}
\author{Hideo Takaoka\thanks{This work was supported by JSPS KAKENHI Grant Number 18H01129.}\\
Department of Mathematics, Kobe University\\
Kobe, 657-8501, Japan\\
takaoka@math.kobe-u.ac.jp}

\date{\empty}

\maketitle

\begin{center}
Dedicated to Professor Tohru Ozawa's on the occasion of his 60th birthday.
\end{center}

\begin{abstract}
In this paper, we consider the bilinear Strichartz estimates for the periodic KdV equation.
We prove the bilinear Strichartz estimate in different frequency regions, which represents a kind of a smoothing effects.
When the period tends to infinity, the frequency localized version of Strichartz estimate on the torus matches to the one in the real line setting.
In addition, the quartic generalized KdV equation is considered as an application.
\end{abstract}

{\it $2020$ Mathematics Subject Classification Numbers.}
35Q55, 42B37.


\section{Introduction}\label{sec:introduction}

This paper is concerned with the periodic Strichartz estimate related to the KdV equation.
Strichartz estimates play an important role in the study on the well-posedness results of dispersive equations.
One of the main examples of dispersive equations is the KdV equation
\begin{eqnarray}\label{eq:kdv}
\partial_tu+\frac{1}{4\pi^2}\partial_x^3u+u\partial_xu=0
\end{eqnarray}
which is a nonlinear evolution model that describes waves in shallow water.
In this paper we discuss the Strichartz for the linear part of the KdV equation, namely, the solution of the initial value problem for the reduced KdV equation
\begin{eqnarray}\label{eq:airy}
\begin{cases}
\partial_tu+\frac{1}{4\pi^2}\partial_x^3u=0,\\
u|_{t=0}=\phi.
\end{cases}
\end{eqnarray}

Let us briefly review the Strichartz estimates of the KdV equation for $\phi\in L^2(\mathbb{R})$, i.e., in the non-periodic case and under initial data decaying at infinity.
The solution of \eqref{eq:airy} is written as
$$
u(t,x)=\left(e^{-\frac{t}{4\pi^2}\partial_x^3}\phi\right)(x)=\int_{\mathbb{R}}e^{2\pi i (\xi^3t+\xi x)}\widehat{\phi}(\xi)\,d\xi
$$
where $\widehat{\phi}$ denotes the Fourier transform of $\phi$. 
In \cite{kpv1}, the following Strichartz estimates are proved by Kenig, Ponce and Vega
\begin{eqnarray}\label{eq:StR} 
\left\||D_x|^{\frac{\vartheta}{4}}e^{-\frac{t}{4\pi^2}\partial_x^3}\phi\right\|_{L_t^q(\mathbb{R}:L^p_x(\mathbb{R}))}\lesssim\|\phi\|_{L^2},
\end{eqnarray}
where
$$
0\le \vartheta\le 1,~q=\frac{4}{\vartheta},~p=\frac{2}{1-\vartheta}.
$$
In the particular case of $\vartheta=2/3$, using the Sobolev inequalities, we obtain the $L^8$ Strichartz estimate
\begin{equation}\label{eq:L8line}
\left\|e^{-\frac{t}{4\pi^2}\partial_x^3}\phi\right\|_{L_{t,x}^8(\mathbb{R}^2)}\lesssim\|\phi\|_{L^2}.
\end{equation}

Similar results to those above hold for the periodic solutions to \eqref{eq:airy}.
Consider the problem of proving the following $L^p$ Strichartz estimates on $\mathbb{T}=\mathbb{R}/\mathbb{Z}$
\begin{eqnarray}\label{eq:st}
\left\|e^{-\frac{t}{4\pi^2}\partial_x^3}\phi\right\|_{L_{t,x}^p(\mathbb{T}^2)}\lesssim\|\phi\|_{L^2},
\end{eqnarray}
where
$$
\left(e^{-\frac{t}{4\pi^2}\partial_x^3}\phi\right)(x)=\sum_{\xi\in\mathbb{Z}}e^{2\pi i(\xi^3 t+\xi x)}\widehat{\phi}(\xi)
$$
and
$$
\widehat{\phi}(\xi)=\int_{\mathbb{T}}e^{-2\pi i\xi x}\phi(x)\,dx.
$$
We begin by summarizing existing results for \eqref{eq:st}.
By the technique used in \cite{z}, we  have that the estimate \eqref{eq:st} holds for $p=4$.
The estimate for $2\le p\le 4$ follows immediately from H\"older's inequality.
In \cite{b1}, Bourgain obtained Strichartz estimates with fractionally loss of derivatives.
Let us recall the general form of those estimates.
We will restrict ourselves to the case of initial data: $\mathrm{supp}\,\widehat{\phi}\subset [-N,N]$ for $N\in\mathbb{N}$.
Then the Strichartz estimates in this direction are the following
\begin{eqnarray}\label{eq:stp}
\left\|\sum_{\xi=-N}^N\widehat{\phi}(\xi)e^{2\pi i(t\xi^3+x\xi)}\right\|_{L^p(\mathbb{T}^2)}\lesssim C(N,p)\left(\sum_{\xi=-N}^N\left|\widehat{\phi}(\xi)\right|^2\right)^{\frac12},
\end{eqnarray}
where $C(N,p)$ is a constant depending only on $p$ and $N$.
It follows from the preceding results that the constant $C(N,p)$ in \eqref{eq:stp} is independent of $N$ if $2\le p\le 4$.
The estimate presented in \cite{b1} is to ensure that the estimate \eqref{eq:stp} holds for $p=6$, up to $\varepsilon$-loss.
Notice that for the non-periodic case, the $L^8$-type Strichartz estimate \eqref{eq:L8line} holds in the absence of dependence on $N$ of $C(N,p)$.
As it was pointed out in \cite{b1}, we expect that the estimate \eqref{eq:stp} holds for $p=8$, up to $C(N,8)=c_{\varepsilon}N^{\varepsilon}$ for any $\varepsilon>0$.
This question is interesting problem, but very hard to answer.

The Strichartz estimate in the range $p\ge 14$  was obtain by Hu and Li \cite{hl}.
They provided the optimal bounds on $C(N,p)$ in \eqref{eq:stp} such that
\begin{eqnarray}\label{eq:Hu}
C(N,p)=c_{\varepsilon}N^{\frac12-\frac{4}{p}+\varepsilon}
\end{eqnarray}
for any $\varepsilon>0$.
The index $1/2-4/p$ of power $N$ is provided by a comparative of $L^8$ Strichartz estimate and of a Sobolev embedding into $L^{\infty}$.
Observe that one has a standard estimate
$$
\left\|\sum_{\xi=-N}^Na_{\xi}e^{2\pi i(t\xi^3+x\xi)}\right\|_{L^{\infty}(\mathbb{T}^2)}\lesssim N^{\frac12+}\left(\sum_{\xi=-N}^N|a_{\xi}|^2\right)^{\frac12},
$$
which is sharp in the sense that the term $N^{1/2+}$ on the right-hand side is optimal.
On the assumption that the estimate \eqref{eq:stp}, up to $\varepsilon$-loss stated above holds, we will expect \eqref{eq:Hu}. 

%

In the present paper, we will discuss the bilinear Strichartz estimates on the large torus.
Let $\mathbb{T}_{\lambda}$ be the torus with the period $\lambda>0$, i.e., $\mathbb{T}_{\lambda}=\mathbb{R}/\lambda\mathbb{Z}$.
Consider the linear $\lambda$ periodic solution to \eqref{eq:airy}
\begin{equation*}
\left(e^{-\frac{t}{4\pi^2}\partial_x^3}\phi\right)(x)=\int e^{2\pi i(\xi^3 t+\xi x)}\widehat{\phi}(\xi)\,\left(d\xi\right)_{\lambda},
\end{equation*}
where
$$
\widehat{\phi}(\xi)=\int_{0}^{\lambda}e^{-2\pi i\xi x}\phi(x)\,dx,\quad \xi\in \mathbb{Z}/\lambda,
$$
and
$$
\int f(\xi)\,\left(d\xi\right)_{\lambda}=\frac{1}{\lambda}\sum_{\xi\in\mathbb{Z}/\lambda}f(\xi).
$$

Let $\eta\in C_0^{\infty}(\mathbb{R})$ be a bump function supported on $[-1,1]$ with $\eta=1$ on $[-1/2,1/2]$.
We may state our result.

\begin{theorem}\label{thm:bilinear}
Let $L$ and $N$ be dyadic numbers satisfying $N\ge L$.
Assume that $\phi_L$ and $\phi_N$ are functions satisfying $\|\phi_L\|_{L^2}=\|\phi_N\|_{L^2}=1$, those Fourier transform have a compact support in $\{\xi\mid |\xi|\sim L\}$ and $\{\xi\mid |\xi|\sim N\}$, respectively.
Then
\begin{equation}\label{eq:bilinear}
\left\| \eta(t)\left(e^{-\frac{t}{4\pi^2}\partial_x^3}\phi_L\right)\left(e^{-\frac{t}{4\pi^2}\partial_x^3}\phi_N\right)\right\|_{L^{4}_{t,x}(\mathbb{R}\times \mathbb{T}_{\lambda})}\lesssim \left(\frac{1}{\lambda}+\frac{1}{N^2}\right)^{\frac{1}{4}}L^{\frac{3}{16}+}N^{\frac{3}{16}+}+\min\left\{\frac{1}{N^{\frac14-}},\frac{1}{\lambda^{\frac14}}+\frac{1}{N^{\frac12}}\right\}L^{\frac14}N^{\frac14}.
\end{equation}
\end{theorem}

\begin{remark}\label{rem:bilinear}
When $N=L$, the bilinear estimate in Theorem \ref{thm:bilinear} is
$$
\left\| \eta(t)e^{-\frac{t}{4\pi^2}\partial_x^3}\phi_N\right\|_{L^{8}_{t,x}(\mathbb{R}\times \mathbb{T}_{\lambda})}\lesssim \left(\frac{1}{\lambda}+\frac{1}{N^2}\right)^{\frac{1}{8}}N^{\frac{3}{16}+}+\min\left\{\frac{1}{N^{\frac18-}},\frac{1}{\lambda^{\frac18}}+\frac{1}{N^{\frac14}}\right\}N^{\frac14}.
$$
In the limit $\lambda\to \infty$, the formal calculations yield
$$
\left\| \eta(t)e^{-\frac{t}{4\pi^2}\partial_x^3}\phi_N\right\|_{L^{8}_{t,x}(\mathbb{R}^2)}\lesssim 1,
$$
which has a close resemblance to the $L^8$ Strichartz estimate \eqref{eq:L8line} in the line case.
\end{remark}

As an application of the above bilinear estimates, we may consider the growth of Sobolev norms of solutions to the quartic generalized KdV equation:
\begin{eqnarray}\label{eq:gkdv}
\begin{cases}
& \partial_t u+\frac{1}{4\pi^2}\partial_x^3 u +{\mathcal P}(u^3)\partial_x u=0, \qquad (t,x)\in\mathbb{R}\times \mathbb{T},\\
& u(0,x) =u_0(x),\qquad x\in\mathbb{T},
\end{cases}
\end{eqnarray}
where ${\mathcal P}$ denotes the orthogonal projection onto mean zero functions
$$
{\mathcal P}(f)(x)=f(x)-\int_{\mathbb{T}}f(x)\,dx.
$$ 
The paper \cite{ckstt} obtained the global well-posedness results as well as the polynomial growth of Sobolev norm of solutions to \eqref{eq:gkdv}, in $H^s(\mathbb{T})$ for $s>5/6$.
This result was improved by Bao and Wu \cite{ba} up to $s\ge 1/2$, where the Sobolev exponent $s=1/2$ is optimal in the sense that the local ill-posedness  holds for $s<1/2$.
One would expect that there is a room for improvement on the growth rate of Sobolev norm of solutions along with the bilinear estimates in Theorem \ref{thm:bilinear}.
In the last section of the paper, we will explore the possibility of adapting the bilinear estimate \eqref{eq:bilinear} to the estimates on the growth order of solutions. 

%

The rest of the paper is organized as follows.
In Section \ref{sec:pre}, we state notation and definitions.
Section \ref{sec:bilinear} proves the bilinear Strichartz estimate in Theorem \ref{thm:bilinear}.
The last section concludes by discussing the growth of Sobolev norms solutions to \eqref{eq:gkdv}. 

\section{Preliminaries}\label{sec:pre}

We list some notation.
The Fourier transform of a function $f$ on $[0,\lambda]$ is defined by
\begin{equation}\label{eq:fourier-s}
\mathcal{F}_{[0,\lambda]}f(\xi)=\widehat{f}(\xi)=\int_{0}^{\lambda}e^{-2\pi i\xi x}f(x)\,dx,\quad \xi\in \mathbb{Z}/\lambda.
\end{equation}
Likewise, we define the Fourier inverse formula defined by
$$
f(x)=\int e^{2\pi i\xi x}\widehat{f}(\xi)\,(d\xi)_{\lambda}=\frac{1}{\lambda}\sum_{\xi\in\mathbb{Z}/\lambda}e^{2\pi i\xi x}\widehat{f}(\xi).
$$
We recall that the following properties:
$$
\|f\|_{L^2[0,\lambda]}=\left\|\widehat{f}\right\|_{L^2((d\xi)_{\lambda})}=\left(\int \left|\widehat{f}(\xi)\right|^2\,(d\xi)_{\lambda}\right)^{1/2},
$$
$$
\int_{[0,\lambda]}f(x)\overline{g(x)}\,dx=\int\widehat{f}(\xi)\overline{\widehat{g}(\xi)}(d\xi)_{\lambda},
$$
$$
\widehat{fg}(\xi)=\int \widehat{f}(\xi-\xi_1)\widehat{g}(\xi_1)\,(d\xi_1)_{\lambda}
$$
and
$$
\widehat{f*g}(\xi)=\widehat{f}(\xi)g(\xi).
$$
We also define the Sobolev space $H^s=H^s(\mathbb{T}_{\lambda})$ with the norm
$$
\|f\|_{H^s}=\left\|\langle \xi\rangle^s\widehat{f}\right\|_{L^2((d\xi)_{\lambda})},
$$
where $\langle \xi\rangle=(1+|\xi|^2)^{1/2}$.

For a function $\phi$ with the time variable $t\in\mathbb{R}$, we define the Fourier transform $\widehat{\phi}$ by
\begin{equation}\label{eq:fourier-t}
\widehat{\phi}(\tau)=\int_{\mathbb{R}}e^{-2\pi i\tau t}\phi(t)\,dt.
\end{equation}
We use the same notation $\widehat{\phi}(\tau)$ here as for the spatial Fourier transform \eqref{eq:fourier-s}, without distinguishing the two.

We use $c,~C$ to denote various independent constants independent of $\lambda$.
We adapt the usual notation that $A\lesssim B$ denotes an estimate of the form $A\le cB$ for some constant $c>0$.
Also we write $A\sim B$ for $A\lesssim B\lesssim A$, and $A=O(B)$ for $|A|\lesssim B$.
The notation $a+$ denotes $a+\varepsilon$ for an arbitrarily small $\epsilon>0$.
Similarly, $a-$ denotes $a-\varepsilon$.

\section{Proof of Theorem \ref{thm:bilinear}}\label{sec:bilinear}

In this section, we prove Theorem \ref{thm:bilinear}.
The proof is due to Bourgain \cite{b1} and also relies on more recent work by Hu and Li \cite{hl}.

\noindent
{\it Proof of Theorem \ref{thm:bilinear}.}
Let $N$ and $L$ be a large integer such that $N\gg L$ and let
$$
\phi_L(x)=\int e^{2\pi i \xi}a_{\xi}\,(d\xi)_{\lambda},\quad 
\phi_N(x)=\int e^{2\pi i \xi}b_{\xi}\,(d\xi)_{\lambda},
$$
where $a_{\xi}=0$ unless $|\xi|\sim L$, and $b_{\xi}=0$ unless $|\xi|\sim N$ satisfying
$$
\int |a_{\xi}|^2\,(d\xi)_{\lambda}=\int |b_{\xi}|^2\,(d\xi)_{\lambda}=1.
$$
It suffices to show that
\begin{equation}\label{eq:suffi}
\|B_{N,L}\|_{L_{t,x}^4}^4\lesssim \left(\frac{1}{\lambda}+\frac{1}{N^2}\right)L^{\frac34+}N^{\frac34+}+\min\left\{\frac{1}{N^{1-}},\frac{1}{\lambda}+\frac{1}{N^2}\right\}LN.
\end{equation}
where
\begin{equation*}
B_{N,L}(t,x)=\eta(t)\int a_{\xi_1}e^{(2\pi i \xi_1)^3 t+2\pi i\xi_1}\,(d\xi_1)_{\lambda}\int \overline{b_{\xi_1}}e^{-(2\pi i\xi_2)^3 t-2\pi i\xi_2}\,(d\xi_2)_{\lambda}.
\end{equation*}

Fix $\mu>0$.
We define
\begin{eqnarray}\label{eq:b3}
f(t,x)=\frac{B_{N,L}(t,x)}{|B_{N,L}(t,x)|}1_{E_{\mu}}(t,x),
\end{eqnarray}
where
$$
E_{\mu}=\{(t,x)\in\mathbb{R}\times\mathbb{T}_{\lambda}\,\mid\, |B_{N,L}(x,t)|\ge \mu\}.
$$
Since
$$
\widehat{B_{N,L}}(\tau,\xi)=\frac{1}{\lambda}\sum_{\xi=\xi_1-\xi_2}\widehat{\eta}\left(\tau+4\pi^2(\xi_1^3-\xi_2^3)\right)a_{\xi_1}b_{\xi_2},
$$
we deduce that
\begin{equation}\label{eq:b4}
\begin{split}
\mu |E_{\mu}|\le & \int_{\mathbb{R}\times \mathbb{T}_{\lambda}}\overline{B_{N,L}(x,t)}f(x,t)\,dxdt\\
= & \frac{1}{\lambda^2}\int_{\mathbb{R}}  \sum_{\scriptstyle \xi\in\mathbb{Z}/\lambda \atop{\scriptstyle \xi=\xi_1-\xi_2}}\overline{\widehat{\eta}(\tau+4\pi^2(\xi_1^3-\xi_2^3))a_{\xi_1}b_{\xi_2}}\widehat{f}(\tau,\xi)\,d\tau. 
\end{split}
\end{equation}
From \eqref{eq:b4}, we obtain
\begin{equation}\label{eq:b6}
\begin{split}
\mu |E_{\mu}| & \lesssim \left(\int |a_{\xi_1}|^2\,(d\xi_1)_{\lambda}\right)^{\frac12}\left(\int |b_{\xi_2}|^2\,(d\xi_2)_{\lambda}\right)^{\frac12}\\
& \times \left(\int\int_{|\xi_1|\sim L}\int_{|\xi_2|\sim N} |\widehat{\eta}(\tau)|\left|\widehat{f}\left(\tau-4\pi^2(\xi_1^3-\xi_2^3),\xi_1-\xi_2\right)\right|^2\,d\tau (d\xi_1)_{\lambda}(d\xi_2)_{\lambda}\right)^{\frac12}.
\end{split}
\end{equation}

We focus our attention on the last term on the right-hand side of \eqref{eq:b6} as
\begin{eqnarray}\label{eq:b7}
\int\int_{|\xi_1|\sim L}\int_{|\xi_2|\sim N} |\widehat{\eta}(\tau)|\left|\widehat{f}(\tau-4\pi^2(\xi_1^3-\xi_2^3),\xi_1-\xi_2)\right|^2\,d\tau (d\xi_1)_{\lambda}(d\xi_2)_{\lambda}.
\end{eqnarray}
Let
$$
\widetilde{\eta}(t)=\int_{\mathbb{R}}e^{2\pi i t\tau}\left|\widehat{\eta}(\tau)\right|\,d\tau
$$
and define the function $K:\mathbb{R}\times \mathbb{T}_{\lambda}\to \mathbb{C}$ by
$$
K(t,x)=\widetilde{\eta}(t)\sum_{\scriptstyle \xi_1,\xi_2\in\mathbb{Z}/\lambda \atop{\scriptstyle |\xi_1|\sim L,|\xi_2|\sim N}}e^{(2\pi i)^3(\xi_1^3-\xi_2^3)t+2\pi i(\xi_1-\xi_2)x}.
$$
The Fourier transformation of $K$ is given by
$$
\widehat{K}(\tau,\xi)=\lambda\sum_{\scriptstyle \xi=\xi_1-\xi_2 \atop{\scriptstyle |\xi_1|\sim L,|\xi_2|\sim N}}\widehat{\widetilde{\eta}}(\tau-4\pi^2(\xi_1^3-\xi_2^3))=\lambda\sum_{\scriptstyle \xi=\xi_1-\xi_2\atop{\scriptstyle |\xi_1|\sim L,|\xi_2|\sim N}}|\widehat{\eta}(\tau-4\pi^2(\xi_1^3-\xi_2^3))|.
$$
From the above, we can thus write \eqref{eq:b7} as
\begin{equation*}
\begin{split}
& \int |\widehat{\eta}(\tau)|\left|\widehat{f}(\tau-4\pi^2(\xi_1^3-\xi_2^3),\xi_1-\xi_2)\right|^2\,d\tau (d\xi_1)_{\lambda}(d\xi_2)_{\lambda}\\
= &\frac{1}{\lambda^2}\int_{\mathbb{R}^2\times (\mathbb{T}_{\lambda})^2}K(t-t',x-x')f(t',x')\overline{f(t,x)}\,dtdt'dxdx'.
\end{split}
\end{equation*}

Let $Q$ be a dyadic number such that
$$
Q=\begin{cases}
N, & \mbox{if $N\lesssim L^2$},\\
L^2, & \mbox{if $N\gg L^2$}.
\end{cases} 
$$
For each $q\sim Q$, we define the set of positive integers less than and relatively prime to $q$ such that
$$
\mathcal{P}_q=\left\{a\in\mathbb{N}\mid 1\le a<q,~\mathrm{gcd}(a,q)=1\right\}.
$$
Recall that the one-dimensional major arcs appear on neighborhoods of sets of rationales
\begin{equation}\label{eq:major}
{\mathcal R}_Q=\left\{\frac{a}{q}~ \middle|~ a\in\mathcal{P}_q,~q\sim Q\right\}.
\end{equation}
Let $\phi$ be a nice bump function supported on $[0,1]$ with $\phi=1$ on $[1/100,2/100]$, and define the function $\Phi:[0,1]\to \mathbb{R}$ to the major arcs by
\begin{equation}
\Phi(t)=\sum_{q\sim Q}\sum_{a\in{\mathcal P}_q}\phi\left(\frac{t-\frac{a}{q}}{\frac{1}{q^2}}\right).
\end{equation}
The bump function $\phi$ is zero everywhere except on $[0,1]$, and this function can be extended to the periodic function.
Upon defining
$$
\mathcal{F}_{[0,1]}\phi(\gamma)=\int_{[0,1]}e^{-2\pi i t\gamma}\phi(t)\,dt,
$$
we have
\begin{equation}\label{eq:FPhi}
\mathcal{F}_{[0,1]}\Phi(\gamma)=\sum_{q\sim Q}\sum_{a\in{\mathcal P}_q}\frac{1}{q^2}e^{-2\pi i\frac{a}{q}\gamma}\mathcal{F}_{[0,1]}\phi\left(\frac{\gamma}{q^2}\right)
\end{equation}
and may write $\Phi$ as Fourier series
$$
\Phi(t)=\sum_{\gamma\in\mathbb{Z}}e^{2\pi it\tau}\mathcal{F}_{[0,1]}\Phi(\gamma).
$$
Recall that
\begin{equation}\label{eq:to}
\sum_{q\sim Q}\sum_{a\in\mathcal{P}_q}\frac{1}{q^2}=\sum_{q\sim Q}\frac{\varphi(q)}{q^2}\sim 1,
\end{equation}
where $\varphi$ is Euler's totient function.
It follows from \eqref{eq:to} that $\mathcal{F}_{[0,1]}\Phi(0)\sim 1$.
Define the following functions
\begin{equation}
\begin{split}
K_1(t,x) & =\frac{\Phi(t)}{\mathcal{F}_{[0,1]}\Phi(0)}K(t,x),\\
K_2(t,x) & = K(t,x)-K_1(t,x).
\end{split}
\end{equation}
Therefore \eqref{eq:b7} becomes the sum of two terms
\begin{equation}\label{eq:K1K2}
\begin{split}
& \frac{1}{\lambda^2}\int_{\mathbb{R}^2\times(\mathbb{T}_{\lambda})^2}K_1(t-t',x-x')f(t',x')\overline{f(t,x)}\,dtdt'dxdx'\\
+& \frac{1}{\lambda^2}\int_{\mathbb{R}^2\times(\mathbb{T}_{\lambda})^2}K_2(t-t',x-x')f(t',x')\overline{f(t,x)}\,dtdt'dxdx'.
\end{split}
\end{equation}
We will give the estimates for \eqref{eq:K1K2} by using interpolation with two upper bound estimates.

Consider the first term in \eqref{eq:K1K2}.
To evaluate $K_1(t,x)$ one need the following Weyl type lemma.
\begin{lemma}[Weyl sum]\label{lem:weyl}
Let $p$ be a natural number, and let $a_2,~a_1,~a_0$ be real numbers.
If $a,~q$ be natural numbers such that $1\le a<q$ with $\mathrm{gcd}(a,q)=1$, then
$$
\left|\sum_{n=0}^p e^{2\pi (tn^3+a_2n^2+a_1n+a_0)}\right|\lesssim p^{1+}\left(\frac{1}{p}+\frac{1}{q}+\frac{q}{p^3}\right)^{\frac14}
$$
for $|t-a/q|\le 1/q^2$.
\end{lemma}
From Lemma \ref{lem:weyl}, we get the pointwise estimate of $K_1(t,x)$.
It follows from the properties for $\mathcal{R}_Q$ in \eqref{eq:major} that the major arcs $a/q+o(1/q^2)$ consist of disjoints intervals in $[0,1]$ for $a/q\in\mathcal{R}_Q$.
Thus from Lemma \ref{lem:weyl}, one has that
\begin{equation}\label{eq:supK1}
\begin{split}
\sup_{(t,x)\in\mathbb{R}\times \mathbb{T}_{\lambda}}|K_1(t,x)|\lesssim  & \lambda^2 L^{1+}N^{1+}\left(\frac{1}{L}+\frac{1}{Q}+\frac{Q}{L^3}\right)^{\frac14}\left(\frac{1}{N}+\frac{1}{Q}+\frac{Q}{N^3}\right)^{\frac14}\\
\lesssim & \lambda^2N^{\frac34+}L^{\frac34+}.
\end{split}
\end{equation}
Using \eqref{eq:supK1}, we have
\begin{equation}\label{eq:first}
\begin{split}
& \frac{1}{\lambda^2}\int_{\mathbb{R}^2\times(\mathbb{T}_{\lambda})^2}K_1(t-t',x-x')f(t',x')\overline{f(t,x)}\,dtdt'dxdx'\\
\lesssim  & N^{\frac34+}L^{\frac34+}\int_{\mathbb{R}^2\times(\mathbb{T}_{\lambda})^2}|f(t',x')f(t,x)|\,dtdt'dxdx\\
\lesssim & N^{\frac34+}L^{\frac34+}|E_{\mu}|^2.
\end{split}
\end{equation}

Considering the second term in \eqref{eq:K1K2}, we need an estimate on
$$
1-\frac{\Phi(t)}{\mathcal{F}_{[0,1]}\Phi(0)}=-\frac{1}{\mathcal{F}_{[0,1]}\Phi(0)}\sum_{\gamma\ne 0}\mathcal{F}_{[0,1]}\Phi(\gamma)e^{2\pi i\gamma t}.
$$
Let $A$ and $B$ be the arbitrary positive constants with $\min\{A,B\}>1$ that we will choose later. 
By \eqref{eq:FPhi}, we have that
\begin{equation*}
\begin{split}
\mathcal{F}_{[0,1]}{\Phi}(\gamma) & \lesssim \sum_{q\sim Q}\left|\sum_{a\in{\mathcal P}_q}\frac{1}{q^2}e^{-2\pi i\left( \frac{a}{q} \right)\gamma}\right|\frac{1}{q^2}\left|\widehat{\phi}\left(\frac{\gamma}{q^2}\right)\right|\\
& \lesssim \frac{1}{Q^2\left\langle\frac{\gamma}{Q^2}\right\rangle^A}\sum_{q\sim Q}\left|\sum_{a\in{\mathcal P}_q}\frac{1}{q^2}e^{-2\pi i\left( \frac{a}{q} \right)\gamma}\right|
\end{split}
\end{equation*}
for $\gamma\ne 0$.
Since by \cite[Lemma 3.33]{b0}, one may bound for the right-hand side of above 
\begin{equation*}
\mathcal{F}_{[0,1]}{\Phi}(\gamma)  \lesssim \frac{d(\gamma,Q)Q^{1+}}{Q^2\left\langle\frac{\gamma}{Q^2}\right\rangle^A}
\lesssim  \frac{\min\{\langle\gamma\rangle^{0+},Q\}}{Q^{1-}\left\langle\frac{\gamma}{Q^2}\right\rangle^A},
\end{equation*}
where $d(\gamma,Q)$ denotes the number of divisors of $\gamma$ less that $Q$.
Since $|\xi_1-\xi_2|\sim |\xi_2|\sim N$ for $|\xi_1|\sim L,~|\xi_2|\sim N$ and $N\gg L$, it follows then from the preceding that
\begin{equation}\label{eq:1b1}
\begin{split}
& \frac{1}{\lambda^2}\sum_{\xi\in\mathbb{Z}/\lambda}\sum_{\scriptstyle \xi=\xi_1-\xi_2 \atop{\scriptstyle |\xi_1|\sim L,|\xi_2|\sim N}}\int_{\mathbb{R}}|\widehat{f}(\tau,\xi)|^2\sum_{0\ne \gamma\in\mathbb{Z}}\mathcal{F}_{[0,1]}\Phi(\gamma)\widehat{\widetilde{\eta}}(\tau-\gamma-4\pi^2(\xi_1^3-\xi_2^3))\,d\tau\\
\lesssim & \frac{1}{\lambda}\|f\|_{L^2}^2\sup_{\scriptstyle (\tau,\xi)\in\mathbb{R}\times \mathbb{Z}/\lambda \atop{\scriptstyle |\xi|\sim N} }\sum_{|\xi_1|\sim L}\sum_{0\ne \gamma\in\mathbb{Z}}\frac{\min\{\langle\gamma\rangle^{0+},Q\}}{Q^{1-}\left\langle\frac{\gamma}{Q^2}\right\rangle^A\langle \tau-\gamma-4\pi^2(\xi_1^3-(\xi_1-\xi)^3)\rangle^B}
\end{split}
\end{equation}
We have to evaluate the contribution of the last term on the right-hand side of \eqref{eq:1b1}. 
Fix $\tau\in\mathbb{R}$ and $\xi\in\mathbb{Z}/\lambda$.
From the identity
$$
\gamma+\left(\tau-\gamma-4\pi^2\left(\xi_1^3-(\xi_1-\xi)^3\right)\right)=\tau-4\pi^2\left(\xi_1^3-(\xi_1-\xi)^3\right),
$$
one has
$$
\max\left\{|\gamma|,\left|\tau-\gamma-4\pi^2\left(\xi_1^3-(\xi_1-\xi)^3)\right)\right|\right\}\gtrsim \left|\tau-4\pi^2\left(\xi_1^3-(\xi_1-\xi)^3\right)\right|,
$$
so that
$$
\sum_{0\ne \gamma\in\mathbb{Z}}\frac{\min\{\langle\gamma\rangle^{0+},Q\}}{Q^{1-}\left\langle\frac{\gamma}{Q^2}\right\rangle^A\langle \tau-\gamma-4\pi^2(\xi_1^3-(\xi_1-\xi)^3)\rangle^B}
\lesssim 
\frac{1}{Q^{1-}\left\langle\frac{\tau-4\pi^2(\xi_1^3-(\xi_1-\xi)^3}{Q^2}\right\rangle^{\min\{A,B\}-}}.
$$
Let $\alpha$ and $\beta$ be roots of the equation $\tau-4\pi^2(\xi_1^3-(\xi_1-\xi)^3)=0$, where $\xi_1$ is a variable, satisfying
$$
\tau-4\pi^2(\xi_1^3-(\xi_1-\xi)^3)=-12\pi^2\xi(\xi_1-\alpha)(\xi_1-\beta).
$$
We may assume here $|\alpha|\ge |\beta|$.
Observe that for $|\xi|\sim N$, 
$$
\langle \tau-4\pi^2(\xi_1^3-(\xi_1-\xi)^3)\rangle\sim \langle N(\xi_1-\alpha)(\xi_1-\beta)\rangle,
$$
and $\alpha,~\beta$ have the form $\alpha+\beta=\xi$, which implies $|\alpha|=\max\{|\alpha|,|\beta|\} \gtrsim |\xi|\sim N$. 
Hence,
\begin{equation*}
\sum_{\scriptstyle \xi_1\in\mathbb{Z}/\lambda \atop{\scriptstyle |\xi_1|\sim L}}\frac{1}{Q^{1-}\left\langle\frac{N(\xi_1-\alpha)(\xi_1-\beta)}{Q^2}\right\rangle^{\min\{A,B\}-}} \lesssim \sum_{\scriptstyle \xi_1\in\mathbb{Z}/\lambda \atop{\scriptstyle |\xi_1|\sim L}}\frac{1}{N^{1-}\langle \xi_1-\beta\rangle} \lesssim \frac{\lambda}{N^{1-}}.
\end{equation*}
Collecting above inequalities, we have that the contribution on the right-hand side of \eqref{eq:1b1} is bounded by
\begin{equation}\label{eq:second}
\begin{split}
& \frac{1}{\lambda^2}\sum_{\xi\in\mathbb{Z}/\lambda}\sum_{\scriptstyle \xi=\xi_1-\xi_2 \atop{\scriptstyle |\xi_1|\sim L,|\xi_2|\sim N}}\int_{\mathbb{R}}|\widehat{f}(\tau,\xi)|^2\sum_{0\ne \gamma\in\mathbb{Z}}\mathcal{F}_{[0,1]}\Phi(\gamma)\widehat{\widetilde{\eta}}(\tau-\gamma-4\pi^2(\xi_1^3-\xi_2^3))\,d\tau\\
\lesssim & \frac{1}{N^{1-}}\|f\|_{L^2}^2\\
= & \frac{1}{N^{1-}}|E_{\mu}|.
\end{split}
\end{equation}
Thus \eqref{eq:first} and \eqref{eq:second} yield
\begin{equation}\label{eq:in1}
\mu^2|E_{\mu}|\lesssim L^{\frac34+}N^{\frac34+}|E_{\mu}|+\frac{1}{N^{1-}}.
\end{equation}

On the other hands, obviously it follows that by straight calculation one has
\begin{equation}\label{eq:in2}
\begin{split}
& \int |\widehat{\eta}(\tau)|\left|\widehat{f}(\tau-4\pi^2(\xi_1^3-\xi_2^3),\xi_1-\xi_2)\right|^2\,d\tau (d\xi_1)_{\lambda}(d\xi_2)_{\lambda}\\
\lesssim & \|f\|_{L^2}^2\frac{1}{\lambda}\sup_{\scriptstyle (\tau,\xi)\in\mathbb{R}\times\mathbb{Z}/\lambda \atop{\scriptstyle |\xi|\sim N}}\sum_{\scriptstyle \xi_1\in\mathbb{Z}/\lambda \atop{\scriptstyle |\xi_1|\sim L}}\frac{1}{\langle \tau+4\pi^2(\xi_1^3-(\xi_1-\xi)^3\rangle^C}\\
\lesssim & \|f\|_{L^2}^2\frac{1}{\lambda}\sup_{ (\tau,\xi)\in\mathbb{R}\times\mathbb{Z}/\lambda }\sum_{\scriptstyle \xi_1\in\mathbb{Z}/\lambda \atop{\scriptstyle |\xi_1|\sim L}}\frac{1}{\langle N^2|\xi_1-\alpha|\rangle^C}\\
\lesssim & \|f\|_{L^2}^2\frac{1}{\lambda}\sum_{\ell\in\mathbb{N}}\frac{1}{\langle \ell\rangle^C}\left\langle \frac{\lambda}{N^2}\right\rangle\\
\lesssim &|E_{\mu}|\left(\frac{1}{\lambda}+\frac{1}{N^2}\right)
\end{split}
\end{equation}
where we choose some $C>1$.

Finally, interpolation between \eqref{eq:in1} and \eqref{eq:in2} yields
\begin{equation}\label{eq:intest}
\mu^2|E_{\mu}|\lesssim \min\left\{L^{\frac34+}N^{\frac34+}|E_{\mu}|+\frac{1}{N^{1-}},\frac{1}{\lambda}+\frac{1}{N^2}\right\}.
\end{equation}

Now come back to the estimate of $\|B_{N,L}\|_{L_{t,x}^4}$ norm.
In the case when $\mu\gg L^{3/8+}N^{3/8+}$, we estimate by \eqref{eq:intest},
$$
\mu^2|E_{\mu}|\lesssim \min\left\{\frac{1}{N^{1-}},\frac{1}{\lambda}+\frac{1}{N^2}\right\}.
$$
Observe that if $\mu\gg L^{1/2}N^{1/2}$, then $E_{\mu}=\emptyset$ since by
$$
\sup_{(t,x)\in\mathbb{R}\times \mathbb{T}_{\lambda}}|B_{N,L}(t,x)|\lesssim L^{\frac12}N^{\frac12}\left(\int |a_{\xi}|^2\,(d\xi)_{\lambda}\right)^{\frac12}\left(\int |b_{\xi}|^2\,(d\xi)_{\lambda}\right)^{\frac12}\le L^{\frac12}N^{\frac12}.
$$
Collecting preceding inequalities, we have
\begin{equation*}
\begin{split}
\|B_{N,L}\|_{L^4_{t,x}}^4 = & \int_0^{CL^{\frac12}N^{\frac12}}\mu^{4-1}|E_{\mu}|\,d\mu\\
\lesssim & \int_{0}^{cL^{\frac38+}N^{\frac38+}}\left(\frac{1}{\lambda}+\frac{1}{N^2}\right)\frac{\mu^{4-1}}{\mu^2}\,d\mu +\min\left\{\frac{1}{N^{1-}},\frac{1}{\lambda}+\frac{1}{N^2}\right\}\int_{cL^{\frac38+}N^{\frac38+}}^{CL^{\frac12}N^{\frac12}}\frac{\mu^{4-1}}{\mu^2}\,d\mu\\
\lesssim & \left(\frac{1}{\lambda}+\frac{1}{N^2}\right)L^{\frac34+}N^{\frac34+}+\min\left\{\frac{1}{N^{1-}},\frac{1}{\lambda}+\frac{1}{N^2}\right\}LN.
\end{split}
\end{equation*}
This completes the proof of Theorem \ref{thm:bilinear}.
\qed

\section{Growth of Sobolev norms of solutions to \eqref{eq:gkdv}}

In this section, we discuss on the growth of Sobolev norms of global in time solutions to the Cauchy problem \eqref{eq:gkdv}.
Our strategy is similar to the one used in the paper of \cite{ckstt}.

Let us review scaling laws.
We rescale the solution by writing
\begin{equation}\label{eq:scal}
u^{\lambda}(t,x)=\frac{1}{\lambda^{\frac23}}u\left(\frac{t}{\lambda^3},\frac{x}{\lambda}\right),
\end{equation}
which maps the function on $\mathbb{R}\times\mathbb{T}$ to that on $\mathbb{R}\times  \mathbb{T}_{\lambda}$.
We shall consider the $\lambda$-periodic Cauchy problem associated to \eqref{eq:gkdv}:
\begin{eqnarray}\label{eq:lgkdv}
\begin{cases}
& \partial_t u^{\lambda}+\frac{1}{4\pi^2}\partial_x^3 u^{\lambda} +{\mathcal P}((u^{\lambda})^3)\partial_x u^{\lambda}=0, \qquad (t,x)\in\mathbb{R}\times \mathbb{T}_{\lambda},\\
& u^{\lambda}(0,x) =u_0^{\lambda}(x),\qquad x\in\mathbb{T}_{\lambda},
\end{cases}
\end{eqnarray}
where we use the same notation ${\mathcal P}$ with no confusion for the orthogonal projection onto mean zero functions on $\mathbb{T}_{\lambda}$ such that
$$
{\mathcal P}(f)(x)=f(x)-\frac{1}{\lambda}\int_{\mathbb{T}_{\lambda}}f(x)\,dx.
$$ 
A straight calculation shows that
$$
\|u^{\lambda}_0\|_{L^2}=\lambda^{-\frac16}\|u_0\|_{L^2}
$$
for $u^{\lambda}_0(x)=\lambda^{-2/3}u_0(x/\lambda)$.

Let $m$ be a smooth non-negative symbol on $\mathbb{R}$, which equals $1$ for $|\xi|< 1$ and equals $|\xi|^{s-1}$ for $|\xi|>2$.
For $N>0$, let $I$ denote the Fourier multiplier operator defined for $f$ belonging to a suitable function class on $\mathbb{T}_{\lambda}$, with symbol such that $\widehat{If}(\xi)=m(\xi/N)\widehat{f}(\xi)$.

The rescaling argument based on \eqref{eq:scal} shows that
\begin{equation}\label{eq:scal1}
\sup_{[0,T]}\|u(t)\|_{H^s}\lesssim \lambda^{\frac16+s}\sup_{[0,\lambda^3T]}\|u^{\lambda}(t)\|_{H^s}\le \lambda^{\frac16+s}\sup_{[0,\lambda^3T]}\|Iu^{\lambda}(t)\|_{H^1}
\end{equation}
and
\begin{equation}\label{eq:scal2}
\|Iu_0^{\lambda}\|_{H^1}\lesssim N^{1-s}\|u_0^{\lambda}\|_{H^s}\lesssim N^{1-s}\lambda^{-\frac16-s}\|u_0\|_{H^s}.
\end{equation}

Let $X_{s,b}$ be the space equipped with the norm
$$
\|u\|_{X_{s,b}}=\left(\int\!\!\int\langle \xi\rangle^{2s}\langle\tau-4\pi^2\xi^3\rangle^{2b}|\widehat{u}(\tau,\xi)|^2\,(d\xi)_{\lambda}d\tau\right)^{\frac12},
$$
where we use the same notation $\widehat{u}(\tau,\xi)$ as for the time-space Fourier transform of a function $u(t,x)$ on $\mathbb{R}\times\mathbb{T}_{\lambda}$, as well.
We also introduce the slightly smaller space $Y^s_{\lambda}$ via the norm
$$
\|u\|_{Y^s_{\lambda}}=\|u\|_{X_{s,\frac12}}+\left(\int\langle\xi\rangle^{2s}\left(\int_{\mathbb{R}}|\widehat{u}(\tau,\xi)|\,d\tau\right)^2\,(d\xi)_{\lambda}\right)^{\frac12}.
$$
The companion space $Z^s$ is defined via the norm
$$
\|u\|_{Z^s_{\lambda}}=\|u\|_{X_{s,-\frac12}}+\left(\int\langle\xi\rangle^{2s}\left(\int_{\mathbb{R}}\frac{|\widehat{u}(\tau,\xi)|}{\langle \tau-4\pi^2\xi^3\rangle}d\tau\right)^2\,(d\xi)_{\lambda}\right)^{\frac12}.
$$
The following lemma follows from Theorem \ref{thm:bilinear}.
\begin{lemma}\label{lem:Fbilinear}
Let $N_1$ and $N_2$ be dyadic numbers satisfying $N_1\ge N_2$.
Assume that $u_{N_1}$ and $u_{N_2}$ are functions on $\mathbb{R}\times \mathbb{T}_{\lambda}$, those Fourier transform have a compact support in $\{\xi\mid |\xi|\sim N_1\}$ and $\{\xi\mid |\xi|\sim N_2\}$, respectively.
Then
\begin{equation}\label{eq:Fbilinear}
\begin{split}
& \left\| u_{N_1}u_{N_2}\right\|_{L^{4}_{t,x}(\mathbb{R}\times \mathbb{T}_{\lambda})}\\
\lesssim& \left[ \left(\frac{1}{\lambda}+\frac{1}{N_1^2}\right)^{\frac{1}{4}}N_2^{\frac{3}{16}+}N_1^{\frac{3}{16}+}+\min\left\{\frac{1}{N_1^{\frac14-}},\frac{1}{\lambda^{\frac14}}+\frac{1}{N_1^{\frac12}}\right\}N_2^{\frac14}N_1^{\frac14}\right]\|u_{N_1}\|_{X_{0,\frac12+}}\|u_{N_2}\|_{X_{0,\frac12+}}.
\end{split}
\end{equation}
\end{lemma}

\noindent
{\it Proof.}
The proof follows by stacking up dyadic block $|\tau-4\pi^2\xi^3|\sim L$ for dyadic number $L$.
\qed

The Cauchy problem \eqref{eq:lgkdv} is locally well-posed in $I^{-1}H^1(\mathbb{T}_{\lambda})$ for $s\ge 1/2$, which was proved in \cite[Lemma 13.1]{ckstt}.

\begin{theorem}[\cite{ckstt}]\label{thm:cksttLWP}
If $s\ge 1/2$, the Cauchy problem \eqref{eq:lgkdv} is locally well-posed for data $u_0^{\lambda}$ satisfying $Iu_0^{\lambda}\in H^1(\mathbb{T}_{\lambda})$.
Moreover, the solution exists on a time interval $[0,\delta]$ with $\delta\sim \|Iu_0^{\lambda}\|_{H^1}^{-\alpha}$ for some $\alpha>0$, and the solution satisfies the estimate
$$
\|Iu^{\lambda}\|_{Y^1_{\lambda}}\lesssim \|Iu_0^{\lambda}\|_{H^1}.
$$
\end{theorem}

Define the Hamiltonian $H(u)$ associated to \eqref{eq:lgkdv} by
\begin{eqnarray}\label{eq:energy}
H(u)=\frac{1}{8\pi^2}\int_{\mathbb{T}_\lambda}(\partial_x u)^2\,dx-\frac{1}{20}\int_{\mathbb{T}_{\lambda}}u^5\,dx,
\end{eqnarray}
which will be a conserved integral under the flow of \eqref{eq:lgkdv}.

We may assume $1/2<s<1$.
Following \cite{ckstt}, we have
\begin{equation}\label{eq:H1data}
H\left(Iu_0^{\lambda}\right)\sim\|I\partial_xu_0^{\lambda}\|_{L^2}^2\sim o(1)
\end{equation}
for choosing
\begin{equation}\label{eq:cN}
N=o\left(\lambda^{(1/6+s)/(1-s)}\right),
\end{equation}
and have the expression
\begin{equation}\label{eq:dhamintonian}
\partial_tH\left(Iu^{\lambda}(t)\right)=\int_{\mathbb{T}_{\lambda}} {\mathcal P}\left[I\left({\mathcal P}((u^{\lambda})^3)\partial_xu^{\lambda}\right)-{\mathcal P}(Iu^{\lambda})^3I\partial_xu^{\lambda}\right)]F \,dx,
\end{equation}
where
$$
F=\frac{1}{4\pi^2}\partial_x^2Iu^{\lambda}+\frac{1}{4}\left(Iu^{\lambda}\right)^4.
$$
The control of the function $F$ of the integrand was proved in \cite[(14.3)]{ckstt}, assuming $\|Iu^{\lambda}\|_{Y^1_{\lambda}}\ll 1$
\begin{equation}\label{eq:ckstt1}
\|F\|_{Y^{-1}_{\lambda}}\ll \lambda^{0+}.
\end{equation}
To evaluate the first term in the integrand, the following estimate was obtained in \cite[(14.5)]{ckstt}.
\begin{lemma}[\cite{ckstt}]\label{kem:ckstt}
Assuming  $\|Iu^{\lambda}\|_{Y^1_{\lambda}}\ll 1$, we have
\begin{equation}\label{eq:ckstt2}
\left\|{\mathcal P}\left[I\left({\mathcal P}((u^{\lambda})^3)\partial_xu^{\lambda}\right)-{\mathcal P}(Iu^{\lambda})^3\partial_xIu^{\lambda}\right)]\right\|_{Z^1_{\lambda}}\ll \lambda^{0+}N^{-\frac12}.
\end{equation}
\end{lemma}
Applying \eqref{eq:H1data}, \eqref{eq:ckstt1} and \eqref{eq:ckstt2} to \eqref{eq:dhamintonian}, one has that from Theorem \ref{thm:cksttLWP}
\begin{equation}\label{eq:ein}
\left|H(Iu^{\lambda}(1))-H(Iu^{\lambda}_0)\right|\ll \lambda^{0+}N^{-\frac12}.
\end{equation}
By iterating the preceding above, we can construct an $H^s(\mathbb{T}_{\lambda})$ solution $u^{\lambda}(t)$ of \eqref{eq:lgkdv} for time $\lambda^{0-}N^{1/2}$.
By \eqref{eq:scal}, we have an $H^s(\mathbb{T})$ solution $u(t)$ for time
$$
T=\lambda^{0-}N^{\frac12}\lambda^{-3}\gtrsim \lambda^{\frac{7s-\frac{35}{6}}{2(1-s)}-}\gg1 
$$
provided $s>5/6$, and the global well-posedness holds on $H^s(\mathbb{T})$ (see \cite[Lemma 13.1]{ckstt}). 

It is worth noting that the exponent $-1/2$ in $N^{-1/2}$ on the right-hand side of \eqref{eq:ckstt2} makes adjustments to an a priori upper bound for the $H^s(\mathbb{T})$ norm of the global in time solutions.
In fact, if the estimate
$$
\left\|{\mathcal P}\left[I\left({\mathcal P}((u^{\lambda})^3)\partial_xu^{\lambda}\right)-{\mathcal P}(Iu^{\lambda})^3\partial_xIu^{\lambda}\right)]\right\|_{Z^1_{\lambda}}\ll \lambda^{0+}N^{-\beta}
$$
holds for some $\beta>0$, then one may have by \eqref{eq:scal1}, \eqref{eq:scal2}, Lemma \ref{thm:cksttLWP}, \eqref{eq:cN} with replacing $\lambda^{3-}T\sim N^{\beta}$ (also $\beta>3(1-s)(1/6+s)$)
$$
\sup_{t\in [0,T]}\|u(t)\|_{H^s(\mathbb{T})}\lesssim T^{\frac{1-s}{\beta-3\frac{1-s}{\frac16+s}}+}.
$$

We shall improve the estimate \eqref{eq:ckstt2} by using the bilinear estimates in Theorem \ref{thm:bilinear}.

\begin{lemma}\label{lem:rbilinear}
For $1/2<s<1$ and the function $u^{\lambda}\in Y^s_{\lambda}$, we have
$$
\left\|{\mathcal P}\left[I\left({\mathcal P}((u^{\lambda})^3)\partial_xu^{\lambda}\right)-{\mathcal P}(Iu^{\lambda})^3\partial_xIu^{\lambda}\right)]\right\|_{Z^1_{\lambda}}\lesssim \lambda^{0+}N^{-1+}\|Iu^{\lambda}\|_{Y^1_{\lambda}}^4.
$$ 
\end{lemma}

\noindent
{\it Proof.}
We repeat the proof of \cite[Lemma 13.1]{ckstt}, but using the bilinear estimates in Theorem \ref{thm:bilinear}.
Following the same scheme, it suffices to show
\begin{equation}\label{eq:suff}
\left\|{\mathcal P}\left[I\left({\mathcal P}(u_1u_2u_3)\partial_xu_4\right)-{\mathcal P}((Iu_1)(Iu_2)(Iu_3))\partial_xIu_4\right)]\right\|_{Z^1_{\lambda}}\lesssim \lambda^{0+}N^{-1+}\prod_{j=1}^4\|Iu_j\|_{Y^1_{\lambda}}.
\end{equation}
Without loss of generality, we may assume that $u_j$ have non-negative Fourier transforms. 
By a Littlewood-Paley decomposition, we restrict each $\widehat{u_j}(\tau_j,\xi_j)$ to Fourier support on the region $|\xi_j|\sim N_j$ for $N_j$ dyadic.  
We prove \eqref{eq:suff} by performing case-by-case analysis.

{\it Case} 1. $N_j\ll N$ for all $1\le j\le 4$: 
By assuming that the functions $u_j$ have Fourier support on $|\xi_j|\ll N$, we easily see that
$$
{\mathcal P}\left[I\left({\mathcal P}(u_1u_2u_3)\partial_xu_4\right)-{\mathcal P}((Iu_1)(Iu_2)(Iu_3))\partial_xIu_4\right]  = {\mathcal P}(u_1u_2u_3)\partial_xu_4-{\mathcal P}(u_1u_2u_3)\partial_xu_4=0.
$$
Then there is no contribution to the left-hand side on \eqref{eq:suff} .

{\it Case} 2. Two of $N_j$ satisfy $N_j\gtrsim N$, and the other two $N_k$ satisfy $N_k\ll N$:   
By the same argument as in the proof of \cite[Lemma 13.1]{ckstt}, it suffices to show that
\begin{equation}\label{eq:suff1}
\left\|{\mathcal P}I\left({\mathcal P}(u_1u_2u_3)\partial_xu_4\right)\right\|_{Z^1_{\lambda}}\lesssim \lambda^{0+}N^{-1+}\prod_{j=1}^4\|Iu_j\|_{Y^1_{\lambda}}.
\end{equation}
Writing
$$
v={\mathcal P}I\left({\mathcal P}(u_1u_2u_3)\partial_xu_4\right),
$$
we have the Fourier transform of $v$ by
$$
\widehat{v}(\tau,\xi)=m(\xi)\int_{\tau=\tau_1+\tau_2+\tau_3+\tau_4}\int_{\scriptstyle \xi=\xi_1+\xi_2+\xi_3+\xi_4 \atop{\scriptstyle \xi_1+\xi_2+\xi_3\ne 0}} i\xi_4\prod_{j=1}^4\widehat{u_j}(\tau_j,\xi_j)\,d\tau_1d\tau_2d\tau_3(d\xi_1)_{\lambda}(d\xi_2)_{\lambda}(d\xi_3)_{\lambda}
$$
for $\xi\ne 0$, and $\widehat{v}(\tau,0)=0$.
The object on right-hand side of \eqref{eq:suff1} is
\begin{equation}\label{eq:suff2}
\left\|\frac{\langle \xi\rangle \widehat{v}(\tau,\xi)}{\langle\tau-4\pi^2\xi^3\rangle^{1/2}}\right\|_{L^2(d\tau(d\xi)_{\lambda})}+\left\|\frac{\langle\xi\rangle\widehat{v}(\tau,\xi)}{\langle\tau-4\pi^2\xi^3\rangle}\right\|_{L^2((d\xi)_{\lambda},L^1(d\tau))}.
\end{equation}
We consider the first term in \eqref{eq:suff2}, because bounding the second term is similar with a small modification of that for the first term as for $s>1/2$. 
By duality, it suffices to show that
\begin{equation}\label{eq:suff3}
\begin{split}
& \int_{\scriptstyle \tau_1+\tau_2+\tau_3+\tau_4+\tau_5=0 \atop{\scriptstyle \xi_1+\xi_2+\xi_3+\xi_4+\xi_5=0\atop{\scriptstyle \xi_5(\xi_1+\xi_2+\xi_3)\ne 0}}} m(\xi_5)\langle\xi_5\rangle|\xi_4|\prod_{j=1}^5\widehat{u_j}(\tau_j,\xi_j)\,d\tau_1d\tau_2d\tau_3d\tau_4(d\xi_1)_{\lambda}(d\xi_2)_{\lambda}(d\xi_3)_{\lambda}(d\xi_4)_{\lambda}\\
 \lesssim & \lambda^{0+}\frac{N_1^{1-}N_2^{1-}N_3^{1-}N_4^{1-}}{N^{1-}}\|Iu_1\|_{Y^0_{\lambda}}\|Iu_2\|_{Y^0_{\lambda}}\|Iu_3\|_{Y^0_{\lambda}}\|Iu_4\|_{Y^0_{\lambda}}\|u_5\|_{X_{0,1/2}}.
\end{split}
\end{equation}
By Leibniz rule, we may assume $N_4\sim N_1\gtrsim N_5\sim N\gg \max\{N_2,N_3\}$, since the other cases are similar.
The inequality \eqref{eq:suff3} is equivalent to
\begin{equation}\label{eq:suff4}
\begin{split}
& \int_{\scriptstyle \tau_1+\tau_2+\tau_3+\tau_4+\tau_5=0 \atop{\scriptstyle \xi_1+\xi_2+\xi_3+\xi_4+\xi_5=0}}  \prod_{j=1}^5\widehat{u_j}(\tau_j,\xi_j)\,d\tau_1d\tau_2d\tau_3d\tau_4(d\xi_1)_{\lambda}(d\xi_2)_{\lambda}(d\xi_3)_{\lambda}(d\xi_4)_{\lambda}\\
\lesssim & \lambda^{0+}\frac{N_1^{2s-1}N_2N_3}{N^{2s-}}\|u_1\|_{Y^0_{\lambda}}\|u_2\|_{Y^0_{\lambda}}\|u_3\|_{Y^0_{\lambda}}\|u_4\|_{Y^0_{\lambda}}\|u_5\|_{X_{0,1/2}}.
\end{split}
\end{equation}
The identities
\begin{equation*}
\sum_{j=1}^5(\tau_j-4\pi^2\xi_j^3)=-12\pi^2\xi_1\xi_4\xi_5+o(N_1^2N),
\end{equation*}
and
\begin{equation*}
\xi_5=-\xi_1-\xi_4+o(N)
\end{equation*}
imply that $N_1\sim N_4$ and
\begin{equation}\label{eq:maxL}
\max_{1\le j\le 5}|\tau_j-4\pi^2\xi_j^3|\gtrsim N_1^2N.
\end{equation}
By another dyadic decomposition, we restrict each $\widehat{u_j}$ to a frequency band  $|\tau_j-4\pi^2\xi_j^3|\sim L_j$ for dyadic $L_j$, and sum in the $L_j$ at the end of the argument.
By \eqref{eq:maxL} we may assume $L_5\gtrsim N_1^2N$, since  the other cases are similar.
Using H\"older inequality, we split the left-hand side of \eqref{eq:suff4} into $L_{t,x}^4,~L^4_{t,x},~L_{x,t}^2$, and apply \eqref{eq:Fbilinear} to control the $L_{t,x}^4$ norms.
From \eqref{eq:cN}, we see $N\ll \lambda\ll N^2$.
Then the contribution of this case to the left-hand side of \eqref{eq:suff4} is bounded by
\begin{equation*}
\begin{split}
\lesssim & \|u_1u_2\|_{L_{t,x}^4}\|u_3u_4\|_{L_{t,x}^4}\|u_5\|_{L_{t,x}^2} \\
\lesssim & \frac{1}{N_1^{1-}N^{\frac12-}}\frac{N_1^{\frac14}N_2^{\frac14}}{\lambda^{\frac14}}\frac{N_3^{\frac14}N_4^{\frac14}}{\lambda^{\frac14}}\prod_{j=1}^5\|u_j\|_{X_{0,1/2-}}\\
\lesssim & \frac{N_1^{2s-1}N_2N_3}{\lambda^{\frac12}N^{2s-}}\prod_{j=1}^5\|u_j\|_{X_{0,1/2}},
\end{split}
\end{equation*}
and the claim \eqref{eq:suff4} follows from $s>1/2$.
%

{\it Case} 3. At least, three of $N_j$ satisfy $N_j\gtrsim N$:
The proof is the same as \cite[Proof of Lemma 13.1]{ckstt}, but there are at least three of $N_j$ satisfying $N_j\gtrsim N$.
By the proof of \cite[Lemma 13.1]{ckstt}, the contribution of this case to the left-hand side of \eqref{eq:suff} is bounded by
\begin{equation}\label{eq:case3}
\lambda^{0+}\|u_1\|_{Y^{1/2}_{\lambda}}\|u_2\|_{Y^{1/2}_{\lambda}}\|u_3\|_{Y^{1/2}_{\lambda}}\|Iu_4\|_{Y^1_{\lambda}}.
\end{equation}
Since at least two of the functions $u_1,~u_2,~u_3$ have the Fourier support in the region $|\xi_j|\gtrsim N$, we use the crude estimate $\|u_j\|_{Y^{1/2}_{\lambda}}\lesssim N^{-1/2}\|Iu_j\|_{Y^1_{\lambda}}$.
Then \eqref{eq:case3} is bounded by
$$
\lesssim \lambda^{0+}N^{-1}\prod_{j=1}^4\|Iu_j\|_{Y^1_{\lambda}},
$$
which is fine.

{\it Case} 4. One of $N_j$ satisfies $N_j\gtrsim N$ and the other three $N_k$ satisfy $N_k\ll N$:
We may assume $N_4\gtrsim N\gg N_3\gtrsim N_2\gtrsim N_1$, since  the other cases are similar.
From the mean value theorem, it follows that
\begin{equation}\label{eq:45-d}
\left| m(\xi_1+\xi_2+\xi_3+\xi_4)-m(\xi_4)\right|\lesssim \frac{|\xi_1+\xi_2+\xi_3|}{N_4}m(N_4).
\end{equation}

{\it Subcase} 4-1: In the region when $N_2\sim N_3$, it follows from \eqref{eq:45-d} that the contribution of this case to the left-hand side of \eqref{eq:suff} is estimated by
\begin{equation*}
\begin{split}
\left\|{\mathcal P}\left({\mathcal P}(u_1u_2u_3)I\langle\partial_x\rangle^{-1}\partial_xu_4\right)\right\|_{Z_{\lambda}^1}
& \sim \left\|{\mathcal P}\left({\mathcal P}(u_1(\langle\partial_x\rangle^{\frac12}u_2)(\langle\partial_x\rangle^{\frac12}u_3))I\langle\partial_x\rangle^{-1/2}\partial_xu_4\right)\right\|_{Z_{\lambda}^{\frac12}}\\
& \lesssim  \lambda^{0+}\|u_1\|_{Y^{\frac12}_{\lambda}}\|\langle\partial_x\rangle^{\frac12}u_2\|_{Y^{\frac12}_{\lambda}}\|\langle\partial_x\rangle^{\frac12}u_3\|_{Y^{\frac12}_{\lambda}}\|I\langle\partial_x\rangle^{-1/2}u_4\|_{Y^{\frac12}_{\lambda}}\\
& \lesssim  \lambda^{0+}N^{-1}\|Iu_1\|_{Y^{1}_{\lambda}}\|Iu_2\|_{Y^{1}_{\lambda}}\|Iu_3\|_{Y^{1}_{\lambda}}\|Iu_4\|_{Y^1_{\lambda}},
\end{split}
\end{equation*}
as before.

{\it Subcase} 4-2: In the region when $N_2\ll N_3$, in a similar way to the Case 2, it suffices to show that
\begin{equation}\label{eq:case4}
\begin{split}
& \int_{\scriptstyle \tau_1+\tau_2+\tau_3+\tau_4+\tau_5=0 \atop{\scriptstyle \xi_1+\xi_2+\xi_3+\xi_4+\xi_5=0\atop{\scriptstyle \xi_5(\xi_1+\xi_2+\xi_3)\ne 0}}} \left(m(\xi_5)-m(\xi_4)\right)\langle\xi_5\rangle|\xi_4|\prod_{j=1}^5\widehat{u_j}(\tau_j,\xi_j)\,d\tau_1d\tau_2d\tau_3d\tau_4(d\xi_1)_{\lambda}(d\xi_2)_{\lambda}(d\xi_3)_{\lambda}(d\xi_4)_{\lambda}\\
 \lesssim & \lambda^{0+}\frac{N_1^{1-}N_2^{1-}N_3^{1-}N_4^{1-}}{N^{1-}}\|Iu_1\|_{Y^0_{\lambda}}\|Iu_2\|_{Y^0_{\lambda}}\|Iu_3\|_{Y^0_{\lambda}}\|Iu_4\|_{Y^0_{\lambda}}\|u_5\|_{X_{0,1/2}}.
\end{split}
\end{equation}
By \eqref{eq:45-d} and $|\xi_1+\xi_2+\xi_3|\sim N_3$, the inequality \eqref{eq:case4} is equivalent to
\begin{equation}\label{eq:case41}
\begin{split}
& \int_{\scriptstyle \tau_1+\tau_2+\tau_3+\tau_4+\tau_5=0 \atop{\scriptstyle \xi_1+\xi_2+\xi_3+\xi_4+\xi_5=0\atop{\scriptstyle \xi_5(\xi_1+\xi_2+\xi_3)\ne 0}}}  \prod_{j=1}^5\widehat{u_j}(\tau_j,\xi_j)\,d\tau_1d\tau_2d\tau_3d\tau_4(d\xi_1)_{\lambda}(d\xi_2)_{\lambda}(d\xi_3)_{\lambda}(d\xi_4)_{\lambda}\\
\lesssim & \lambda^{0+}\frac{N_1N_2}{N^{1-}N_4^{0+}}\|u_1\|_{Y^0_{\lambda}}\|u_2\|_{Y^0_{\lambda}}\|u_3\|_{Y^0_{\lambda}}\|u_4\|_{Y^0_{\lambda}}\|u_5\|_{X_{0,1/2}}.
\end{split}
\end{equation}
The identity
\begin{equation*}
\sum_{j=1}^5(\tau_j-4\pi^2\xi_j^3)=12\pi^2(\xi_1+\xi_2+\xi_3)(\xi_4^2+\xi_5^2+\xi_4\xi_5)+O(N_3^3),
\end{equation*}
implies that
\begin{equation*}
\max_{1\le j\le 5}|\tau_j-4\pi^2\xi_j^3|\gtrsim N_4^2N_3.
\end{equation*}
By using dyadic decomposition, we again restrict each $\widehat{u_j}$ to a frequency band  $|\tau_j-4\pi^2\xi_j^3|\sim L_j$ for dyadic $L_j$, and sum in the $L_j$ at the end of the argument.
We may only consider the region when $L_5=\max\{L_1,L_2,L_3,L_4,L_5\}\gtrsim N_4^2N_3$, as before.
Using H\"older inequality, we split the left-hand side of \eqref{eq:case41} into $L_{t,x}^{\infty},~L^{\infty}_{t,x},~L_{x,t}^{\infty},~L_{t,x}^2,~L_{t,x}^2$.
Then the contribution of this case to the left-hand side of \eqref{eq:case41} is bounded by
\begin{equation*}
\begin{split}
\lesssim & \|u_1\|_{L_{t,x}^{\infty}} \|u_2\|_{L_{t,x}^{\infty}} \|u_3\|_{L_{t,x}^{\infty}} \|u_4\|_{L_{t,x}^{2}} \|u_5\|_{L_{t,x}^{2}}\\
\lesssim & \frac{N_1^{\frac12+}N_2^{\frac12+}N_3^{\frac12+}}{N_4^{1-}N_3^{\frac12-}}\prod_{j=1}^5\|u_j\|_{X_{0,1/2-}}\\
\lesssim & \frac{N_1N_2}{N^{1-}N_4^{0+}}\prod_{j=1}^5\|u_j\|_{X_{0,1/2}},
\end{split}
\end{equation*}
which is fine.

This completes the proof of Lemma \ref{lem:rbilinear}.
\qed

\begin{remark}
The paper \cite{ba} improved the a priori estimate \eqref{eq:ein} by adding appropriate correction term to $H(Iu)$.
The advantage of the bilinear estimate \eqref{eq:bilinear} is that we may improve the a priori estimate obtained in \cite{ba}.
We will discuss this problem in a forthcoming paper.
\end{remark}

\bigskip

\noindent
{\it Acknowledgments.}
The author would like to thank Professor Tadahiro Oh for providing me some related papers.
The prior version of  this paper was written to explore the counterexample of $L^p$ Strichartz estimates for $p=8$, but our analyses were methodologically flawed, the conclusions were incorrect.
Accordingly, the paper was revised. 
The author would like to thank the anonymous reviewers for many helpful comments and suggestions to the prior version of  this paper .
%
%


%

\end{document}